\documentclass[a4paper,10pt]{article}
\usepackage{amsmath}
\usepackage{amssymb}
\usepackage{graphics}
\usepackage{rotating}
\usepackage{cite}
\usepackage{color}

\title{A generalization of a Cullen's Integral Theorem for the quaternions}
\author{Daniel Alay\'{o}n-Solarz \underline{(danieldaniel@gmail.com)}}

\begin{document}

\maketitle

\begin{abstract}
We discuss the proof of a certain integral theorem obtained by C. G. Cullen, originally stated on the class of the analytic intrinsic functions on the quaternions. It is shown that this integral theorem is true for a larger class of quaternionic functions. 
\end{abstract}

\section{Introduction and preliminaries}
Let \textbf{H} be the algebra of the quaternions and let $p$ be a quaternion, then $p$ can be written as
\begin{displaymath}
p := t + xi + yj + zk.
\end{displaymath}
Let $f$ be a quaternion-valued function of a single quaternionic variable. Consider the class of complex-like quaternionic functions such that a member $f$ can be written as 
\begin{displaymath}
f = u + \iota v,
\end{displaymath}
where $u$ and $v$ are real functions and $\iota$ is defined as
\begin{displaymath}
\iota := \frac{xi+yj+zk}{\sqrt{x^{2}+y^{2}+z^{2}}}.
\end{displaymath}
In particular the identity function, sending one quaternion onto itself is complex-like and thus the quaternion $p$ can be written as:
\begin{displaymath}
p = t + r \iota,
\end{displaymath}
where
\begin{displaymath}
r:=\sqrt{x^{2}+y^{2}+z^{2}}.
\end{displaymath}
Recall the left-Fueter operator is given by: 
\begin{displaymath}
 D_{l}:=\frac{\partial  }{\partial t}+ i\frac{\partial }{\partial x} + j\frac{\partial }{\partial y} + k\frac{\partial }{\partial z}. 
\end{displaymath}
For the rest of this paper, the class of complex-like quaternionic functions satisfying
\begin{displaymath}
 D_{l}f=\frac{-2v}{r}
\end{displaymath}
which will be referred as \textbf{Hyperholomorphic}. Note that $\iota$ can be parametrized by spherical coordinates
\begin{displaymath}
\iota = (\cos\alpha \sin\beta, \sin\alpha \sin\beta, \cos \beta).
\end{displaymath}
The coordinate system based in the variables $(t,r,\alpha,\beta)$ is especially well suited to study the interplay between the Fueter operator and the complex-like quaternionic functions. The Fueter operator in this coordinate system has the form:
\begin{displaymath}
D_{l} = \frac{\partial }{\partial t} + \iota \frac{\partial }{\partial r}  - \frac{1}{r} \frac{\partial}{\partial_{l} \iota},
\end{displaymath}
where
\begin{equation*}
\frac{\partial}{\partial_{l} \iota} := ({\iota}_{\alpha})^{-1}\frac{\partial}{\partial \alpha} +  ({\iota}_{\beta})^{-1}\frac{\partial}{\partial \beta}.
\end{equation*}
and $\iota_{\alpha}$ and $\iota_{\alpha}$ represent the derivatives of $\iota$ with respect to $\alpha$ and $\beta$ respectively.
Using  this coordinate system the following characterization of Hyperholomorphic functions can be proved:
\newtheorem{teo2}{Proposition}
\begin{teo2}
A function $f=u + \iota v$ is hyperholomorphic if and only if $u$ and $v$ satisfy:
\begin{eqnarray*}
\frac{\partial u}{\partial t} - \frac{\partial v}{\partial r} = 0,
\end{eqnarray*}
\begin{eqnarray*}
\frac{\partial v}{\partial t} +\frac{\partial u}{\partial r}= 0,
\end{eqnarray*}
and
\begin{eqnarray*}
\frac{\partial v}{\partial \alpha}(\sin\beta)^{-1}+\frac{\partial u}{\partial \beta} = 0,
\\
\frac{\partial u}{\partial \alpha}(\sin\beta)^{-1}-\frac{\partial v}{\partial \beta}= 0.
\end{eqnarray*}
\end{teo2}
It is well known that regular functions in the Fueter sense, that is, quaternionic null-solutions to the Fueter operator, are in general not closed under the quaternionic product. However non-zero hyperholomorphic functions form a multiplicative group:

\begin{teo2}
The sum and product of two hyperholomorphic functions is hyperholomorphic. If a function is hyperholomorphic and non-zero then its algebraic inverse is hyperholomorphic. 
\end{teo2}
A important property of hyperholomorphic functions is that they generalize the Fueter's Theorem \cite{F35}:
\begin{teo2}
Let $f$ be a hyperholomorphic function. Then
\begin{equation*}
D_{l}\Delta{f}=D_{r}\Delta{f}=0.
\end{equation*}
\end{teo2} 
Where by $D_{r}$ we denote the right-Fueter operator. Fueter's theorem was originally stated for the smaller class of quaternionic functions obtained by rotating around the real axis a complex analytic function. \\
The purpose of this paper is to show how hyperholomorphic functions satisfy the following integral theorem given by Cullen in \cite{Cul}.
\begin{teo2}
Let $f= u + \iota v$ be an hyperholomorphic function and let $K$ any smooth, simple closed hypersurface in H the quaternionic space, disjoint from the real axis, $K^{*}$ being the interior of $K$. Let $n(p)=n_{0}+n_{1}i+n_{2}j+n_{3}k$ where $(n_{0},n_{1},n_{2},n_{3})$ is the unit outer normal to $K$ at $p$. Then,
\begin{displaymath}
\int_{K}n(p) f(p)\frac{1}{r^2}dS_{K}= -2 \int_{K^{*}} u \frac{\iota}{r^3} dV.
\end{displaymath}
where $dS_{K}$ is the element of surface area on $K$.
\end{teo2}

\section{The Proof}
This proof we will consider is the same Cullen gave in \cite{Cul} and it requires two preliminar lemmas. Our contribution consists on the observation that one of these lemmas used by Cullen and satisfied by the class of analytic intrinsic functions is actually a definition of the larger class of hyperholomorphic functions. 
\begin{teo2}(Cullen's lemma)
A function $f = u + \iota v$ a hyperholomorphic  if and only if:
\begin{displaymath}
D_{l}(\frac{f}{r^2}) = -\frac{2}{r^{3}}u \iota.
\end{displaymath}
\end{teo2}
\textbf{Proof.}
We first assume $f$ is hyperholomorphic, then:
\begin{displaymath}
D_{l}(\frac{f}{r^2})= (\frac{\partial }{\partial t} + \iota \frac{\partial }{\partial r}  - \frac{1}{r} \frac{\partial}{\partial_{l} \iota}) (\frac{u+ \iota v}{r^2})
\end{displaymath}
\begin{displaymath}
= (\frac{\partial }{\partial t} + \iota \frac{\partial }{\partial r})(\frac{f}{r^2}) - \frac{2}{r^3}v = \frac{1}{r^2}(\frac{\partial f}{\partial t} + \iota \frac{\partial f}{\partial r}) - \frac{2}{r^3}(u \iota -v) - \frac{2}{r^3}v=  -\frac{2}{r^{3}}u \iota.
\end{displaymath}
Note how we use the fact that since $f$ is hyperholomorphic then:
\begin{displaymath}
\frac{\partial f}{\partial_{l} \iota}  = 2v
\end{displaymath}
which in turn is a consequence of the hyperholomorphic functions being  Cullen-regular  \cite{GG}, that is it satisfies the following equation:
\begin{displaymath}
( \frac{\partial }{\partial t} + \iota \frac{\partial }{\partial r} )f = 0
\end{displaymath}
Note that for all functions $f=u+\iota v$:
\begin{displaymath}
D_{l}(\frac{f}{r^2}) = \frac{1}{r^2}( \frac{\partial }{\partial t} + \iota \frac{\partial }{\partial r}  - \frac{1}{r} \frac{\partial}{\partial_{l} \iota})f - \frac{2}{r^3}(\iota f)
\end{displaymath}
so in particular if $f$ satisfies the Cullen lemma then:
\begin{displaymath}
\frac{1}{r^2}( \frac{\partial }{\partial t} + \iota \frac{\partial }{\partial r}  - \frac{1}{r} \frac{\partial}{\partial_{l} \iota})f - \frac{2}{r^3}(\iota f) = -\frac{2}{r^{3}}u \iota
\end{displaymath}
implies
\begin{displaymath}
\frac{1}{r^2}D_{l}f = -\frac{2}{r^3}v
\end{displaymath}
which for $r \neq 0$ is the condition for hyperholomorphicity.
We continue with the second lemma:
\begin{teo2}
Let $(n_{0},n_{1},n_{2},n_{3})$ be the outward unit normal of $K$, then
\begin{displaymath}
\int_{K^{*}}(\frac{\partial f_{0}}{\partial t}+\frac{\partial f_{1}}{\partial x}+\frac{\partial f_{2}}{\partial y}+\frac{\partial f_{3}}{\partial z}) dV = \int_{K}(f_{0}n_{0}+f_{1}n_{1}+f_{2}n_{2}+f_{3}n_{3}) dS_{k}
\end{displaymath}
where $f_{i}$ are differentiable quaternionic functions.
\end{teo2}
\textbf{Proof.} Following Cullen, this is an application of the Gauss Theorem in four dimensions for the components of $f_{i}$. \\
We are now ready to prove our main result, for this it suffices to show that for $f=u+\iota v$ an hyperholomorphic function we have:
\begin{displaymath}
\int_{K} n(p)f(p)\frac{1}{r^2}dS_{K} = \int_{K}(n_{0}\frac{f}{r^2}+n_{1}\frac{if}{r^2}+n_{2}\frac{jf}{r^2}+n_{3}\frac{kf}{r^2}) dS_{K} = \int_{K^{*}} D_{l
} (\frac{f}{r^2}) dV
\end{displaymath}
\begin{displaymath}
 = -2 \int_{K^{*}} u \frac{\iota}{r^3} dV.
\end{displaymath} 
And it is proved. \\
If one starts with the class of right-hyperholomorphic functions, that is those functions of the form $f=u+\iota v$ that satisfy
\begin{displaymath}
 D_{r}f=\frac{-2v}{r}
\end{displaymath}
then the integral theorem would read:
\begin{displaymath}
\int_{K}f(p)n(p)\frac{1}{r^2}dS_{K}= -2 \int_{K^{*}} u \frac{\iota}{r^3} dV.
\end{displaymath}
The class of functions that are both left- and right-hyperholomorphic is not empty. For example the class of analytic instrinsic functions on the quaternions is of this form.





\end{document}